\documentclass[11pt,dvips]{amsart}  
\usepackage{amssymb,amsmath,amsfonts}  
\usepackage[dvips]{graphicx,xcolor}
\usepackage{cite}
\allowdisplaybreaks[4]
\definecolor{pgray}{gray}{0.8}

\newtheorem{theorem}{Theorem}[section]

\newtheorem{remark}[theorem]{Remark}
\numberwithin{equation}{section}

\title{\mbox{}}%Asymptotic behaviour of solutions to a parabolic system of semi-conductor model}
\begin{document}
\begin{center}
{\bf \LARGE{
	Distortion of charge distribution due to internal electric fields described by the drift-diffusion semiconductor model\footnote{In this manuscript, some carelessness are remained. The corrected version is available in Z. Angew. Math. Phys.  \textbf{77} (2026), 148.
	\texttt{https://doi.org/10.1007/s00033-026-02800-1}}
}}\\
\vspace{5mm}
Masakazu Yamamoto\footnote{e-mail : \texttt{mk-yamamoto@gunma-u.ac.jp}}\\Graduate School of Science and Technology, Gunma University
\end{center}
\maketitle
\vspace{-15mm}
%%% ----------------------------------------------------------------------
%
%% main text
\begin{abstract}
In this paper, the initial value problem for the Debye--H\"uckel drift-diffusion equation is studied.
This equation was introduced as a model describing plasma behavior and is also known as a simulation model of MOSFET, and so its solution describes charge density.
It is well-known that, if the initial density is localized, then the density is adjusted to be radially symmetric due to the linear diffusion.
Consequently, the electric field is also governed by a radially symmetric potential, and its effects are expected to act radially symmetrically.
The main result express the electric field and its effect on the charge density as concrete functions.
It also denotes the distortion of symmetry and the shift of scale on the density due to the internal electric field.
Unlike the historical paper via Escobedo and Zuazua and the followers, the main result captures stronger nonlinearity than the logarithmic shift.
%Particularly, this effect also is radially symmetric.
%This result is consistent with the fact that a radially symmetric distribution of charge forms a radially symmetric electric field.
\end{abstract}

\section{Introduction}
We study large-time behavior of solutions of the following initial-value problem of the drift-diffuion equation of Debye--H\"uckel and Nernst--Planck type:
\begin{equation}\label{dd}
\left\{
\begin{array}{lr}
	\partial_t u - \Delta u = \nabla \cdot (u\nabla\psi),
	&
	t > 0,~ x \in \mathbb{R}^3,\\
	-\Delta\psi = u,
	&
	t > 0,~ x \in \mathbb{R}^3,\\
	u(0,x) = u_0 (x),
	&
	x \in \mathbb{R}^3,
\end{array}
\right.
\end{equation}
where $(u,\psi) = (u,\psi) (t,x)$ stand for density of charge and potential of internal electric field, respectively, and $u_0 = u_0 (x)$ is a given initial density.
The drift-diffusion equation was developed as a model of plasma dynamics.
This equation also models the diffusion of ions in solution.
A significant application of this equation is the simulation model for MOSFETs (cf.\cite{Fng-Ito,Jngl,Mck}).
From this perspective, it is more appropriate to consider bipolar types rather than monopolar types.
In fact, for large-time behavior of solutions, it has been proven that there is no difference in their mathematical treatment except in electrically neutral cases (see\cite{KwshmKbysh}).
Well-posedness of \eqref{dd} is confirmed in \cite{KO}.
For positive solutions, Moser--Nash energy theory guarantees the global in time solvability and the decay estimates that
\begin{equation}\label{decay}
	\| u(t) \|_{L^q (\mathbb{R}^3)}
	\le
	C (1+t)^{-\gamma_q}
\end{equation}
for $1 \le q \le \infty$ and $\gamma_q = \frac32 (1-\frac1q)$.
Moreover
\[
	\| u(t) - M_0 G (t) \|_{L^q (\mathbb{R}^3)}
	= O (t^{-\gamma_q-\frac12}\log t)
\]
as $t \to +\infty$ for $M_0 = \int_{\mathbb{R}^3} u_0 (x) dx$ and Gaussian $G(t,x) = (4\pi t)^{-3/2} e^{-|x|^2/(4t)}$ (cf.\cite{BlrDlbl,Jngl}).
This estimate states that the solution is approximated asymptotically by the solution of heat equation.
In other words, the charge distribution is corrected to be radially symmetric.
The indicator developed by Escobedo and Zuazua\cite{EZ} is used to predict whether this estimate is optimal.
According to this indicator, it is expected that there is some constant $d \in \mathbb{R}^3$ such that
\[
	u(t) \sim M_0 G(t) + d \cdot \nabla G(t) \log t
\]
as $t \to +\infty$.
If we put $K(t) = d \cdot \nabla G(t) \log t$, then $K$ represents scale shift in the solution caused by nonlinear effects.
Moreover, this satisfies $\partial_t K = \Delta K - \frac{K}{t\log t}$.
Therefore, if this expectation is true, then the solution is approximated by the linear solutions.
In fact, this expectation is false.
Using the renormalization developed in \cite{KtM,Ymd}, Ogawa and the author\cite{OgY} showed that $d = 0$ and derived the higher-order expansion.
Namely,
\begin{equation}\label{OgYexp}
	u(t) \sim M_0 G(t) + M_1 \cdot \nabla G(t) + M_0^2 J(t)
\end{equation}
as $t \to +\infty$ for $M_0 = \int_{\mathbb{R}^3} u_0 (x) dx,~ M_1 = - \int_{\mathbb{R}^3} x u_0 (x) dx$ and
\[
	J(t) = \int_0^t \nabla G(t-s) * (G\nabla (-\Delta)^{-1} G) (1+s) ds.
\]
At first glance, $J$ appears to be logarithmically growing since
\[
	\| J(t) \|_{L^q (\mathbb{R}^3)} \le Ct^{-\gamma_q -\frac12} \log (1+t).
\]
However, since $\int_{\mathbb{R}^3} (G\nabla (-\Delta)^{-1}G) dy = 0$, this function is represented by
\[
	J(t) = \int_0^t \int_{\mathbb{R}^3} (\nabla G(t-s,x-y) - \nabla G(t-s,x)) \cdot (G\nabla (-\Delta)^{-1}G) (1+s,y) dyds,
\]
and the mean value theorem deletes the logarithm.
Namely,
\begin{equation}\label{estJ}
	\| J(t) \|_{L^q (\mathbb{R}^3)} = O (t^{-\gamma_q-\frac12})
\end{equation}
and then
\begin{equation}\label{exp1st}
	\| u(t) - M_0 G(t) \|_{L^q (\mathbb{R}^3)} = O (t^{-\gamma_q-\frac12})
\end{equation}
as $t \to +\infty$.
The expansion \eqref{OgYexp} provides a precise approximation of the solution, but it leaves some ambiguity in the spatiotemporal structure.
The term $J$ describes the nonlinear effects on the solution.
However, its large-time behavior is not as clear as seen in \eqref{estJ}.
Furthermore, its spatial structure is complex.
As seen in \eqref{exp1st}, if the charge density is nearly radially symmetric, the electric field should also act radially symmetrically.
Indeed, we establish our main result as follows.
\begin{theorem}\label{thm}
Let $u_0 \in L^1 (\mathbb{R}^3) \cap L^\infty (\mathbb{R}^3)$ be positive and $x u_0 \in L^1 (\mathbb{R}^3)$.
Then
\[
	\biggl\| u(t) - M_0 G(t) - M_1 \cdot \nabla G(t) - \frac{M_0^2}{8\pi^2} \int_0^t \int_0^\infty s^{-1/2} (2+\sigma)^{-5/2} \Delta G \left( t - \tfrac{s}{2+\sigma} \right) d\sigma ds \biggr\|_{L^q (\mathbb{R}^3)}
	= o (t^{-\gamma_q-\frac12})
\]
as $t\to +\infty$ for $M_0 = \int_{\mathbb{R}^3} u_0 (x) dx$ and $M_1 = - \int_{\mathbb{R}^3} x u_0 (x) dx$, and $1 \le q \le \infty$ and $\gamma_q = \frac32 (1-\frac1q)$.
In addition, if $|x|^2 u_0 \in L^1 (\mathbb{R}^3)$, then
\[
\begin{split}
	&\biggl\| u(t) - M_0 G(t) - M_1 \cdot \nabla G(t) - \frac{M_0^2}{8\pi^2} \int_0^t \int_0^\infty s^{-1/2} (2+\sigma)^{-5/2} \Delta G \left( t - \tfrac{s}{2+\sigma} \right) d\sigma ds\\
	&\hspace{15mm} + \frac{(2\pi-3\sqrt{3}) M_0^3}{2^{7}\cdot 3^2\cdot \pi^3} \Delta G(t) \log t \biggr\|_{L^q (\mathbb{R}^3)}
	= O (t^{-\gamma_q-1})
\end{split}
\]
as $t\to +\infty$.
\end{theorem}
In the expansion,
\[
	U_0 (t) = M_0 G(t)\quad\text{and}\quad
%\]
%and
%\[
	U_1^{\mathrm{odd}} (t)
	= M_1 \cdot \nabla G(t)
\]
describe linear diffusions, and
\begin{equation}\label{defU1rad}
	U_1^{\mathrm{rad}} (t)
	=
	\frac{M_0^2}{8\pi^2} \int_0^t \int_0^\infty s^{-1/2} (2+\sigma)^{-5/2} \Delta G \left( t - \tfrac{s}{2+\sigma} \right) d\sigma ds
\end{equation}
and
\[
	K_2 (t) \log t = -\frac{(2\pi-3\sqrt{3}) M_0^3}{2^{7}\cdot 3^2\cdot \pi^3} \Delta G(t) \log t
\]
capture the characteristics of internal electric field, respectively.
Firstly, $U_1^{\mathrm{odd}}$ and $U_1^{\mathrm{rad}}$ have the same scale that is
\[
	\lambda^{3+1} (U_1^{\mathrm{odd}}, U_1^{\mathrm{rad}}) (\lambda^2 t, \lambda x) = (U_1^{\mathrm{odd}}, U_1^{\mathrm{rad}}) (t, x)
\]
for $\lambda > 0$.
This property immediately gives that
\[
	t^{\gamma_q+\frac12} \| (U_1^{\mathrm{odd}},U_1^{\mathrm{rad}}) (t) \|_{L^q (\mathbb{R}^3)} =  \| (U_1^{\mathrm{odd}},U_1^{\mathrm{rad}}) (1) \|_{L^q (\mathbb{R}^3)}
\]
for $t > 0$ and $1 \le q \le \infty$, where $\gamma_q = \frac32 (1-\frac1q)$.
However, their spatial symmetries are entirely different.
This fact violates the second law of thermodynamics and the principle of entropy increase.
In general linear diffusion phenomena, the decay rate of solution in time clearly corresponds to its spatial symmetry.
More precisely, $U_1^{\mathrm{rad}}$ does not decay as rapidly as might be expected given its high symmetry.
Expanding $\Delta G (t-\frac{s}{2+\sigma})$ around $-\frac{s}{2+\sigma} = 0$ rewrites $U_1^{\mathrm{rad}}$ as
\[
	U_1^{\mathrm{rad}} (t) = - \frac{M_0^2}{4\pi^2} \sum_{n=0}^\infty \frac1{(2n+3)(2n+1)} \left( \frac{t}2 \right)^{n+\frac12}\frac{(-\Delta)^{n+1} G(t)}{n!}
\]
and all terms decay with the same rate.
On the other hand, the logarithmic shift $K_2 \log t$ has no linear counterpart for comparison, making it an isolated nonlinearity.
From these, we understand that these functions cannot arise in linear phenomena and $U_1^{\mathrm{rad}}$ represents the largest component among the nonlinear effects.
These terms represent symmetry distortions and scale shifts in the charge density due to internal electric field.
%From another perspective, $U_1 = U_1^{\mathrm{odd}} + U_1^{\mathrm{rad}}$ is corresponds to the Cole--Hopf solution in the Burgers equation.
Secondary, not only the largest nonlinear component $U_1^{\mathrm{rad}}$ but also the second component $K_2\log t$ is radially symmetric.
This indicates that the correction \eqref{exp1st} to the charge density strongly influences the electric field.
Since $\lambda^{3+2} K_2 (\lambda^2 t, \lambda x) = K_2 (t,x)$ for $\lambda > 0$, this term has a weaker influence than $U_1^{\mathrm{rad}}$.
We will see that this weak shift $K_2 \log t$ is generated from the strong distortion $U_1^{\mathrm{rad}}$.
To explain the above discussion in semiconductor engineering terms, $U_0$ and $U_1^{\mathrm{odd}}$ represent linear diffusion effects, while $U_1^{\mathrm{rad}}$ and $K_2 \log t$ describe field effects.
Furthermore, $U_1^{\mathrm{odd}}$ represents uneven distribution of charge, and the influence of the electric field generates a radially symmetric component $U_1^{\mathrm{rad}}$ that balances it.
We also emphasize that these terms are determined only by the moments of the initial-data.
Especially, the field effects are fixed only by the mass of charge.
This is a characteristic not observed in convection diffusion phenomena.
Our problem has the same scale as
\[
	\partial_t u - \Delta u = d \cdot \nabla (|u|^{2/3}u)
\]
in $\mathbb{R}^3$, where $d$ is some direction.
Compare our results, for example, with the assertions in \cite{DrCrpo,EZ,Ksb}.

Theorem \ref{thm} treats positive solutions.
In fact, the same assertion fulfills even for solutions where the sign changes if the initial-data is small.
This theorem actually requires \eqref{decay}.
For the same reason, the similar expansion is derived for small solutions to the Keller--Segel chemotaxis model and one to the gravitational interaction model in astronomy.
These models are formulated as $\partial_t u - \Delta u + \nabla \cdot (u\nabla\psi) = 0$ for $-\Delta \psi = u$.

\vspace{2mm}

\paragraph{\textbf{Notations.}}
We often omit the spatial parameter from functions, for example, $u(t) = u(t,x)$.
In particular, $G(t) * u_0 = \int_{\mathbb{R}^3} G(t,x-y) u_0 (y) dy$ and $\int_0^t g(t-s) * f(s) ds = \int_0^t \int_{\mathbb{R}^3} g(t-s,x-y) f(s,y) dyds$.
We symbolize the derivations by $\partial_t = \partial/\partial t,~ \partial_j = \partial/\partial x_j$ for $j = 1,2,3,~ \nabla = (\partial_1,\partial_2,\partial_3)$ and $\Delta = \partial_1^2 + \partial_2^2 + \partial_3^2$.
The length of a multiindex $\alpha = (\alpha_1,\alpha_2, \alpha_3) \in \mathbb{Z}_+^3$ is given by $\lvert \alpha \rvert = \alpha_1 + \alpha_2 + \alpha_3$, where $\mathbb{Z}_+ = \mathbb{N} \cup \{ 0 \}$.
We abbreviate that $\alpha ! = \alpha_1 ! \alpha_2! \alpha_3 !,~ x^\alpha = x_1^{\alpha_1} x_2^{\alpha_2} x_3^{\alpha_3}$ and $\nabla^\alpha = \partial_1^{\alpha_1} \partial_2^{\alpha_2}  \partial_3^{\alpha_3}$.
We define the Fourier transform and its inverse by $\mathcal{F} [\varphi] (\xi) = (2\pi)^{-3/2} \int_{\mathbb{R}^3} \varphi (x) e^{-ix\cdot\xi} dx$ and $\mathcal{F}^{-1} [\varphi] (x) = (2\pi)^{-3/2} \int_{\mathbb{R}^3} \varphi (\xi) e^{ix\cdot\xi} d\xi$, respectively, where $i = \sqrt{-1}$.
Therefore, the field generated by a charge density $g$ is given by $\nabla (-\Delta)^{-1} g = \mathcal{F}^{-1} [i\xi |\xi|^{-2} \hat{g}]$.
The Lebesgue space and its norm are denoted by $L^q (\mathbb{R}^3)$ and $\| \cdot \|_{L^q (\mathbb{R}^3)}$, that is, $\| f \|_{L^q (\mathbb{R}^3)} = (\int_{\mathbb{R}^3} |f(x)|^q dx)^{1/q}$ for $1 \le q < \infty$ and $\| f \|_{L^\infty (\mathbb{R}^3)}$ is the essential supremum.
The heat kernel and its decay rate on $L^q (\mathbb{R}^3)$ are symbolized by $G(t,x) = (4\pi t)^{-3/2} e^{-|x|^2/(4t)}$ and $\gamma_q = \frac{3}2 (1-\frac1q)$.
%In our definitions, $G(t) = (2\pi)^{-3/2} \mathcal{F}^{-1} [e^{-t|\xi|^2}]$.
%We denote the floor function by the Gauss symbol $[\mu] = \max \{ m \in \mathbb{Z} \mid m \le \mu \}$.
We employ Landau symbol.
Namely, $f(t) = o(t^{-\mu})$ and $g(t) = O(t^{-\mu})$ mean $t^\mu f(t) \to 0$ and $t^\mu g(t) \to c$ for some $c \in \mathbb{R}$ such as $t \to +\infty$ or $t \to +0$, respectively.
A subscript of function represents its scale or decay rate.
For example, $\| U_m (t) \|_{L^q (\mathbb{R}^n)} = t^{-\gamma_q-\frac{m}2} \| U_m (1) \|_{L^q (\mathbb{R}^n)}$ for $t > 0$, and $\| r_m (t) \|_{L^q (\mathbb{R}^n)} = O (t^{-\gamma_q-\frac{m}2})$ or $O(t^{-\gamma_q-\frac{m}2} \log t)$ as $t \to +\infty$.
Various positive constants are simply denoted by $C$.

\section{Proof of main result}
Our proof is based on Escobedo--Zuazua theory with the renormalization for the mild solution
\begin{equation}\label{ms}
	u(t) = G(t) *u_0 + \int_0^t \nabla G(t-s) * (u\nabla (-\Delta)^{-1} u) (s)ds.
\end{equation}
As well known, the term of initial-data is expanded as $G(t) * u_0 \sim U_0 (t) + U_1^{\mathrm{odd}} (t)$ as $t \to +\infty$.
Since $\int_{\mathbb{R}^3} (u\nabla (-\Delta)^{-1}u) dy = 0$ and $\int_{\mathbb{R}^3} (G\nabla (-\Delta)^{-1} G) dy = 0$, the nonlinear term is expanded as
\[
	\int_0^t \nabla G(t-s) * (u\nabla (-\Delta)^{-1} u) (s)ds
	=
	U_1^{\mathrm{rad}} (t) + r_2 (t)
\]
for
\begin{equation}\label{repU1rad}
%\[
	U_1^{\mathrm{rad}} (t)
	=
	M_0^2 \int_0^t \nabla G(t-s) * (G\nabla (-\Delta)^{-1} G) (s) ds
%\]
\end{equation}
of the original form and
\[
	r_2 (t) =
	\int_0^t \int_{\mathbb{R}^3}
		(\nabla G(t-s,x-y) - \nabla G(t-s,x))
		\cdot (u\nabla (-\Delta)^{-1} u - U_0 \nabla (-\Delta)^{-1} U_0) (s,y)
	dyds.
\]
The same procedure as in \cite{OgY} gives that
\begin{equation}\label{estr2}
	\| r_2 (t) \|_{L^q (\mathbb{R}^3)} = O (t^{-\gamma_q-1} \log t)
\end{equation}
as $t\to +\infty$.
Consequently, we see the first assertion if the representation \eqref{defU1rad} is guaranteed.
From
\[
	\nabla (-\Delta)^{-1} G(t)
	= \frac{2\sqrt{\pi}}\pi\int_0^\infty \nabla G(t+\sigma) d\sigma,
\]
we see that
\[
\begin{split}
	&(G\nabla(-\Delta)^{-1}G) (t)
	=
%	(4\pi)^{-3/2} t \int_0^\infty (2t+\sigma)^{-5/2} \nabla G \left( \tfrac{t(t+\sigma)}{2t+\sigma} \right) d\sigma\\
%	&=
	(4\pi)^{-3/2} t^{-1/2} \int_0^\infty (2+\sigma)^{-5/2} \nabla G \left( \tfrac{t(1+\sigma)}{2+\sigma} \right) d\sigma.
\end{split}
\]
This rewrites \eqref{repU1rad} to \eqref{defU1rad}.
We derive the logarithmic evolution and prove the second assertion.
The remaind term is further expanded as
\[
\begin{split}
	&r_2 (t) = 
	- \sum_{|\beta|=1} \nabla^\beta \nabla G(t) \cdot \int_0^t \int_{\mathbb{R}^3} y^\beta (U_0\nabla(-\Delta)^{-1}U_1 + U_1 \nabla (-\Delta)^{-1} U_0) (1+s,y) dyds
	 + \tilde{r}_2 (t),
\end{split}
\]
where $U_1 = U_1^{\mathrm{odd}} + U_1^{\mathrm{rad}}$ and
\[
\begin{split}
	 &\tilde{r}_2 (t)
	=
	- \sum_{|\beta|=1} \nabla^\beta \nabla G(t) \cdot \int_0^t \int_{\mathbb{R}^3} y^\beta \bigl( (u\nabla(-\Delta)^{-1}u - U_0 \nabla (-\Delta)^{-1} U_0) (s,y)\\
	 &\hspace{15mm}  -  (U_0\nabla(-\Delta)^{-1}U_1 + U_1 \nabla (-\Delta)^{-1} U_0) (1+s,y) \bigr) dyds\\
	 &+ \int_0^t \int_{\mathbb{R}^3} \biggl( \nabla G(t-s,x-y) - \sum_{|\beta|=0}^1 \nabla^\beta \nabla G(t,x) (-y)^\beta \biggr)\cdot (u\nabla(-\Delta)^{-1}u - U_0 \nabla (-\Delta)^{-1} U_0) (s,y) dyds.
\end{split}
\]
Focusing on these symmetries reveals that several terms in the first part of $r_2$ vanish.
Indeed, we see that
\[
\begin{split}
	&\int_0^t \nabla G(t-s) * (u\nabla (-\Delta)^{-1} u) (s)ds
	=
	U_1^{\mathrm{rad}} (t)\\
	&- \frac13 \Delta G(t) \int_0^t (1+s)^{-1} ds \int_{\mathbb{R}^3} y \cdot (U_0\nabla(-\Delta)^{-1}U_1^{\mathrm{rad}} + U_1^{\mathrm{rad}} \nabla (-\Delta)^{-1} U_0) (1,y) dy
	 + \tilde{r}_2 (t).
\end{split}
\]
Here we separated the integration by using the scales of $U_0$ and $U_1^{\mathrm{rad}}$.
We calculate the coefficient, then
\[
\begin{split}
	 &\int_{\mathbb{R}^3} y \cdot (U_0\nabla(-\Delta)^{-1}U_1^{\mathrm{rad}} + U_1^{\mathrm{rad}} \nabla (-\Delta)^{-1} U_0) (1,y) dy\\
	 &= -\frac{\sqrt{\pi}}{8\pi^2} M_0^3 \int_0^1 \int_0^\infty s^{-1/2} (2+\sigma)^{-5/2} \int_{\mathbb{R}^3} G(1,y) y \cdot \nabla G \left( 1 - \tfrac{s}{2+\sigma},y \right) dyd\sigma ds\\
	 &+\frac{\sqrt{\pi}}{8\pi^2} M_0^3 \int_0^1 \int_0^\infty s^{-1/2} (2+\sigma)^{-5/2}  \int_{\mathbb{R}^3} \Delta G \left( 1 - \tfrac{s}{2+\sigma},y \right) y \cdot \nabla (-\Delta)^{-1} G(1,y)  dyd\sigma ds.
\end{split}
\]
Hence, the harmonic analysis and the elementary calculus say that
\[
\begin{split}
	 &\int_{\mathbb{R}^3} y \cdot (U_0\nabla(-\Delta)^{-1}U_1^{\mathrm{rad}} + U_1^{\mathrm{rad}} \nabla (-\Delta)^{-1} U_0) (1,y) dy\\
	 &=
	 \frac{M_0^3}{64\pi^3} \int_0^1 \int_0^\infty s^{-1/2} (2+\sigma)^{-1} (4+2\sigma -s)^{-3/2} d\sigma ds
	 =
	 \frac{(2\pi -3\sqrt{3})M_0^3}{2^{7}\cdot 3 \cdot \pi^3},
\end{split}
\]
and we deside the coefficient of logarithmic shift.
On the same way as in \cite{OgY}, we have that
\begin{equation}\label{esttr2}
	\| \tilde{r}_2 (t) \|_{L^q (\mathbb{R}^3)} = O (t^{-\gamma_q-1})
\end{equation}
as $t \to +\infty$.
We review this procedure for reader's convenience.
Upon the condition $|x|^2 u_0 \in L^1 (\mathbb{R}^3)$, it is well known that
\[
	\| G(t) * u_0 - M_0 G(t) - M_1 \cdot \nabla G(t) \|_{L^q (\mathbb{R}^3)} = O (t^{-\gamma_q-1})
\]
as $t \to +\infty$.
A coupling of this and \eqref{estr2} yields that
\[
	\| u(t) - U_0 (t) - U_1 (t) \|_{L^q (\mathbb{R}^3)} = O (t^{-\gamma_q-1}\log t)
\]
as $t \to +\infty$.
This estimate together with Hardy--Littlewood--Sobolev inequality that
\[
	\| \nabla (-\Delta)^{-1} \varphi \|_{L^{q_*} (\mathbb{R}^3)} \le C \| \varphi \|_{L^q (\mathbb{R}^3)}
\]
for $1 < q< 3$ and $\frac1{q_*} = \frac1q - \frac13$, and the weighted estimates that
\begin{equation}\label{decay-wt}
	\| |x|^m u(t) \|_{L^q (\mathbb{R}^3)} = O (t^{-\gamma_q+\frac{m}2})
\end{equation}
as $t \to +\infty$  yield \eqref{esttr2}.
Indeed, the coefficient of first term of $\tilde{r}_2$ fulfills
\[
\begin{split}
	&\biggl| \int_0^t \int_{\mathbb{R}^3} y^\beta \bigl( (u\nabla(-\Delta)^{-1}u - U_0 \nabla (-\Delta)^{-1} U_0) (s,y)\\
	 &\hspace{15mm}  -  (U_0\nabla(-\Delta)^{-1}U_1 + U_1 \nabla (-\Delta)^{-1} U_0) (1+s,y) \bigr) dyds \biggr|\\
	 &\le
	 \int_0^t \int_{\mathbb{R}^3} \bigl( \bigl| y^\beta u(s,y) \bigr| + \bigl| y^\beta U_0 (s,y) \bigr| \bigr) \bigl| \nabla (-\Delta)^{-1} \bigl( (u-U_0) (s,y) - U_1 (1+s,y) \bigr) \bigr| dyds\\
	 &+  \int_0^t \int_{\mathbb{R}^3} \bigl| y^\beta U_1 (1+s,y) \bigr| \bigl| \nabla (-\Delta)^{-1} (u-U_0) (s,y) \bigr| dyds
	\le
	C \int_0^t s^{-1/2} (1+s)^{-1} \log (2+s) ds \le C.
\end{split}
\]
The details of Hardy--Littlewood--Sobolev inequality are found in\cite{Stin,Zmr}.
For the proof of the weighted estimates, see \cite{OgY}.
Taylor theorem rewrites the second term of $\tilde{r}_2$ as
\[
\begin{split}
	&\int_0^t \int_{\mathbb{R}^3} \biggl( \nabla G(t-s,x-y) - \sum_{|\beta|=0}^1 \nabla^\beta \nabla G(t,x) (-y)^\beta \biggr)\cdot (u\nabla(-\Delta)^{-1}u - U_0 \nabla (-\Delta)^{-1} U_0) (s,y) dyds\\
	&= \sum_{|\beta|=2} \int_0^{t/2} \int_{\mathbb{R}^3} \int_0^1 \frac{\nabla^\beta \nabla G(t- s, x-\lambda y)}{\beta!} \lambda (-y)^\beta \cdot (u\nabla(-\Delta)^{-1}u - U_0 \nabla (-\Delta)^{-1} U_0) (s,y) d\lambda dyds\\
	&+ \sum_{|\beta|=0}^1 \int_0^{t/2} \int_{\mathbb{R}^3} \int_0^1 \partial_t \nabla^\beta \nabla G(t-\lambda s,x) (-s) \cdot (u\nabla(-\Delta)^{-1}u - U_0 \nabla (-\Delta)^{-1} U_0) (s,y) d\lambda dyds\\
	&+ \int_{t/2}^t \int_{\mathbb{R}^3} \biggl( \nabla G(t-s,x-y) - \sum_{|\beta|=0}^1 \nabla^\beta \nabla G(t,x) (-y)^\beta \biggr)\cdot (u\nabla(-\Delta)^{-1}u - U_0 \nabla (-\Delta)^{-1} U_0) (s,y) dyds.
\end{split}
\]
A combination of $L^p$-$L^q$ estimate, \eqref{exp1st}, \eqref{decay-wt} and Hardy--Littlewood--Sobolev inequality gives that
\[
\begin{split}
	&\biggl\| \sum_{|\beta|=2} \int_0^{t/2} \int_{\mathbb{R}^3} \int_0^1 \frac{\nabla^\beta \nabla G(t- s, x-\lambda y)}{\beta!} \lambda (-y)^\beta \cdot (u\nabla(-\Delta)^{-1}u - U_0 \nabla (-\Delta)^{-1} U_0) (s,y) d\lambda dyds \biggr\|_{L^q (\mathbb{R}^3)}\\
	&\le
	C \int_0^{t/2} (t-s)^{-\gamma_q - \frac32} s^{-\frac12} ds
	= Ct^{-\gamma_q-1}.
\end{split}
\]
The next part is treated on the same way, and
\[
\begin{split}
	&\biggl\| \int_{t/2}^t \int_{\mathbb{R}^3} \biggl( \nabla G(t-s,x-y) - \sum_{|\beta|=0}^1 \nabla^\beta \nabla G(t,x) (-y)^\beta \biggr)\\
		&\hspace{15mm} \cdot (u\nabla(-\Delta)^{-1}u - U_0 \nabla (-\Delta)^{-1} U_0) (s,y) dyds \biggr\|_{L^q (\mathbb{R}^3)}\\
	&\le
	C \int_{t/2}^t (t-s)^{-\frac12} s^{-\gamma_q - \frac32} ds + C\sum_{|\beta|=0}^1 t^{-\gamma_q-\frac{|\beta|}2-\frac{1}2} \int_{t/2}^t s^{-\frac32+\frac{|\beta|}2} ds
	= Ct^{-\gamma_q-1}.
\end{split}
\]
Therefore we get \eqref{esttr2} and complete the proof.

\begin{remark}
Our assertions and \eqref{OgYexp} are not contradictory.
Indeed, from \eqref{repU1rad}, we see that
\[
\begin{split}
%	&\rho_2 (t) = 
	&M_0^2 J(t) - U_1^{\mathrm{rad}} (t)\\
	&= \int_0^t \int_{\mathbb{R}^3} (\nabla G(t-s,x-y) - \nabla G(t-s,x))
	 \cdot \left( (G\nabla (-\Delta)^{-1}G) (1+s,y) - (G\nabla (-\Delta)^{-1}G) (s,y) \right) dyds
\end{split}
\]
since $\int_{\mathbb{R}^3} (G\nabla (-\Delta)^{-1}G) dy = 0$.
Hence, the $L^p$-$L^q$ estimate together with the mean value theorem gives that $\| M_0^2 J(t) - U_1^{\mathrm{rad}} (t) \|_{L^q (\mathbb{R}^3)} = O (t^{-\gamma_q-1} \log t)$ as $t \to +\infty$.
The first assertion also can be shown from this estimate.
\end{remark}
\begin{remark}
Needless to say, $K_2 \neq 0$ as long as $M_0 \neq 0$.
Hence, tracing the origin of $K_2$ reveals that $U_1^{\mathrm{rad}} \neq 0$.
\end{remark}

\end{document}